\documentclass{article}

\begin{document}

\title{A note on the $\alpha$-invariant of the Mukai-Umemura 3-fold }
\author{S. K. Donaldson}
\maketitle
\newtheorem{thm}{Theorem}
\newtheorem{lem}{Lemma}
\newtheorem{cor}{Corollary}
\newtheorem{prop}{Proposition}
\newcommand{\bC}{{\bf C}}
\newcommand{\bR}{{\bf R}}
\newcommand{\cH}{{\cal H}}
\newcommand{\cP}{{\cal P}}
\newcommand{\bP}{{\bf P}}
\newcommand{\oa}{\overline{a}}
\newcommand{\ua}{\underline{a}}

In this note we discuss a $3$-dimensional Fano $3$-fold $X$ discovered by
Mukai and Umemura \cite{kn:MU}. The relevant properties of $X$, for us, are the following
\begin{itemize}
\item The anti-canonical bundle $K_{X}^{-1}$ is very ample.
\item There is a holomorphic action of $PSL(2,\bC)$ on $X$ and a point $x_{0}\in
X$ with stabiliser the iscosahedral group $\Gamma\subset SO(3)\subset PSL(2,\bC)$.
\item There is a $PSL(2,\bC)$-invariant section $s$ of $K_{X}^{-1}$ whose
zero-set $\Sigma \subset X$ is the complement of the $PSL(2,\bC)$-orbit of
$x_{0}$.
\item The section $s$ vanishes with multiplicity $1$ on $\Sigma$ and $\Sigma$
has at worst cusp-like singularities in that, around any non-smooth point of $\Sigma$,
we can find local complex  co-ordinates $z_{1}, z_{2}, z_{3}$ such that
$\Sigma$ is defined by the equation $z_{1}^{2}= z_{2}^{3}$.
\end{itemize}
These are all well-known facts: for completeness we give a brief account
in Section 3.

\

For any Fano $n$-fold $Z$, with a compact group $G$ of holomorphic automorphisms,
Tian defined an invariant $\alpha_{G}(Z)$\cite{kn:T1}. This is a positive  real number
and Tian showed that if $\alpha_{G}(Z)>\frac{n}{n+1}$ then $Z$ has a Kahler-Einstein
metric. In our case, we 
fix the compact group $SO(3)\subset PSL(2,\bC)$ containing the stabiliser of $x_{0}$. We prove
\begin{thm}
The $\alpha$-invariant $\alpha_{SO(3)}(X)$ is $5/6$.
\end{thm}
We should say straightaway that this is not really a new result. Alessio
Corti has explained to the author that, given the facts above, it can be obtained from the more general
theories of \cite{kn:DKol}. But the extreme simplicity (almost triviality) of the argument 
we give,
for this special case, made the author feel that it might be worth writing
it down. Our approach extends, in an obvious way, to cover  other situations,
if suitable examples can be found. In Section 3 we explain how the same argument
gives another proof of the result of Batyrev and Selinova, for symmetric toric manifolds.

\

Of course the point of Theorem 1 is that, since $5/6> 3/4$, Tian's work gives
\begin{cor} The Mukai-Umemura manifold $X$ has a Kahler-Einstein metric.
\end{cor}
This is a particularly interesting example since, as Tian shows in \cite{kn:T2},
there are arbitrarily small deformations of $X$ which do not admit Kahler-Einstein
metrics. The situation is analogous to that occurring for bundles over curves
(say). If $E_{0}$ is the trivial  vector bundle  over a curve, of rank $2$
or more, then $E_{0}$ admits a flat unitary connection but, if the genus
of the curve is not zero, there are arbitrarily small deformations of $E_{0}$
which do not admit flat unitary connections. (See also the discussion of
this phenomenon in
\cite{kn:DK}, 6.4.3, Example (ii)).

\

The author is grateful to Alessio Corti, Julien Keller and Frances Kirwan
for  helpful advice and  discussions.

\

\section{Proof of the Theorem}

We will only write down the proof that $\alpha\geq 5/6$, which is what is
relevant to Corollary 1. The proof that $\alpha=5/6$ is an easy extension
of this.

\

We begin by recalling the definition of Tian's $\alpha$-invariant.
Let $Z$ be a Fano manifold on which a compact group $G$ acts by holomorphic
automorphisms and fix a $G$-invariant Kahler metric $\omega_{0}$ in the
cohomology class $-c_{1}(K_{Z})$. Let $\cP$ be the set of $G$-invariant Kahler
potentials $\phi$ on $X$ such that $\omega_{\phi}= \omega_{0}+i \partial
\overline{\partial} \phi>0$ and $\max_{Z} \phi=0$. Thus $\cP$ can be identified
with
the set of all $G$-invariant Kahler metrics in the given Kahler class. Let
$A\subset \bR$ be the set defined by the condition that $\beta\in A$ if there
exists a $C_{\beta}\in \bR$ such that
$$ \int_{Z} e^{-\beta \phi} d\mu_{0} \leq C_{\beta}, $$
for all $\phi\in \cP$. Here $d\mu_{0}$ is the volume form defined by the fixed
metric $\omega_{0}$. Then Tian sets
$$  \alpha_{G}(Z) = \sup\{\beta: \beta \in A\}, $$
and shows that this does not depend on the choice of $\omega_{0}$.

\

It will be important to have the correct signs in our discussion, so we note
that $$\Lambda i\partial \overline{\partial}\phi = \frac{1}{2} \Delta_{0} \phi$$ where $$
\Lambda:\Omega^{1,1}\rightarrow \Omega^{0}$$ is the inner product with the
metric form $\omega_{0}$ and  we use the sign convention that, on Euclidean space, the Laplacian is $\Delta \phi =\sum \frac{\partial^{2}}{\partial x_{i}^{2}}$.

\

\begin{lem}
There is an $M\in \bR$ such that
$$\int_{Z} \phi \ d\mu_{0} \geq -M$$
for all $\phi\in \cP$.
\end{lem}
This is a step in Tian's proof that $\alpha>0$ and we repeat his argument.
If $\phi\in \cP$ we have
$$  \Delta_{0} \phi = 2\Lambda( i\partial\overline{\partial} \phi)\geq
-2n. $$
Let $K$ be the Green's function for $\Delta_{0}$, so that for all functions
$f$ on $Z$
$$  f(x) = -\int_{Z} K(x,y) (\Delta_{0} f)(y) d\mu_{0}(y)+\frac{1}{V} \int
f(y) d\mu_{0}(y), $$
where $V$ is the volume of the manifold. With our sign conventions, $K$ is bounded below and, since we can change $K$ by the addition of a constant without
affecting the identity, we may suppose that $K\geq 0$. While $K$ is singular
along the diagonal it is  integrable in each variable. Let $x$ be the point
where $\phi$ vanishes. Then applying the Green's identity to $\phi$ we have
$$\int_{Z} \phi(y) d\mu_{0}(y) = V \int_{Z} K(x,y) \Delta_{0}\phi d\mu_{0}(y)\geq
-2n V \int_{Z} K(x,y) d\mu_{0}(y). $$
So we can take $$M= 2n V \max_{x}\int_{Z} K(x,y) d\mu_{0}(y). $$

\

\

For the rest of this section we work with  the Mukai-Umemura manifold
$X$. There is a Hermitian metric  on the line bundle $K_{X}^{-1}$ such
that the curvature of the associated unitary connection is $-2\pi i \omega_{0}$.
Set $$f_{0}=- \log \left(\vert s \vert^{2}\right).$$
This is a smooth function on $X\setminus \Sigma$ and $i\partial \overline{\partial}
f_{0} = 2\pi \omega_{0}$.
 
\begin{lem}
For any $\beta<\frac{5}{6}$ the function $\exp(\beta f_{0})$ is integrable.
\end{lem}

This is also standard. The integral in question is
$$  \int_{Z} \vert s\vert^{-2\beta} d\mu_{0}. $$
By our assumptions we can reduce to considering the
integrals
$$  \int_{B^{4}} \vert z^{2} - w^{3} \vert^{-2\beta} , $$
where $B^{4}$ is the unit ball in $\bC^{2}$ and $z,w$ are complex co-ordinates.
 For integer $r\geq 0$ let $\Omega_{r}$ be the annular region
$$ \Omega_{r} = \{ (z,w) \in \bC^{2}: 2^{-3(r+1)} \leq\vert z\vert\leq 2 ^{-3r}, 2^{-2(r+1)}\leq\vert w\vert\leq 2^{-2r}\}, $$
and let
$$I_{r} = \int_{\Omega_{r}}  \vert z^{2} - w^{3} \vert^{-2\beta} . $$
The substitution $z'= 2^{3} z, w'= 2^{2} w$ shows that
$$  I_{r+1} = 2^{(12 \beta- 10)} I_{r}. $$
Thus $\sum_{r} I_{r}$ is finite if $\beta<5/6$ and the union of the $\Omega_{r}$
cover $B^{4}\setminus \{0\}$.

\

Now we give the main proof. Let $\phi$ be an $SO(3)$-invariant function in $\cP$.
Then $f= f_{0}+2\pi \phi$ is an $SO(3)$-invariant function on $X\setminus \Sigma$
with $i\partial \overline{\partial} f =2\pi\omega_{\phi}>0$. Define a function
$\tilde{f}$ on $PSL(2,\bC)$ by the obvious procedure, i.e.
$$  \tilde{f}(g)= f(g(x_{0})), $$
for $g\in PSL(2,\bC)$.
Then if $u$ lies in $SO(3)$ and $\gamma$ lies in $\Gamma$ it follows from
the definitions that $\tilde{f}(u g \gamma) =\tilde{f}(g)$. Thus
$\tilde{f}$ can be regarded as a function on the homogeneous space
$\cH= PSL(2,\bC)/SO(3)$ which is invariant under the natural action of $\Gamma$
on $\cH$. Write $P_{0}$ for the base-point in $\cH$ corresponding to the
identity in $PSL(2,\bC)$. Then $P_{0}$ is the unique fixed point for the
$\Gamma$-action on $\cH$. Of course $\cH$ is just a model of real hyperbolic $3$- space, and $\Gamma$ acts on $\cH$ by rotations about $P_{0}$.  
Now we have
\begin{prop}
The function $\tilde{f}$  on $\cH$ is convex along geodesics.
\end{prop}
This is a very standard fact. Indeed, any time we have a compact Lie group
$G$ with complexification $G^{c}$ we can the identify  $G$-invariant plurisubharmonic functions
on $G^{c}$ with convex functions on $G^{c}/G$. The operator $i\partial \overline{\partial}$
on $G^{c}$, when applied to a $G$-invariant function, can be identified with the Hessian operator on $G^{c}/G$. We omit the easy calculation.

\

Now our function $\tilde{f}$ on $\cH$ is proper and bounded below, by construction, so it achieves
a minimum in $\cH$. By  the convexity and $\Gamma$-invariance this minimum
must occur at $P_{0}$. (This observation makes up the main content of our
argument.)  Set $\tilde{f}(P_{0})=-2\pi b$. Then the inequality $\tilde{f}\geq
-2\pi b$ translates back into the statement that $\phi\geq f_{0} -b$. So
$$  \int_{Z} e^{-\beta \phi} d\mu_{0}\leq e^{b \beta} \int_{Z} f_{0}^{-\beta}d\mu_{0}.
$$
By Lemma 1, it suffices to obtain an upper bound on $b=b(\phi)$ for all $\phi\in \cP$. Let $\tilde{f}_{0}$ be the function on $\cH$ defined in the same way as $\tilde{f}$
but using $f_{0}$. Thus $\tilde{f} \leq \tilde{f_{0}}$, since $\phi\leq 0$.
Let $D$ be the geodesic ball in $\cH$ centred on $P_{0}$, of radius $1$ say. Let $2\pi \oa$ be the maximum value of $\tilde{f}_{0}$ on $D$, so for any $\phi$
we have $\tilde{f}\leq 2\pi \oa$ on $D$. Convexity along geodesics emanating from
$P_{0}$ implies that
$$  \frac{1}{2\pi}\tilde{f}(Q)\leq -b+ (\oa+b)\ {\rm dist}(Q,P_{0}), $$
for any point $Q$ in $D$. In particular, on the ball $\frac{1}{2} D$ of radius $1/2$ about $P_{0}$ we have $\tilde{f} \leq \pi(\oa-b)$. 

Take the inverse image in $PSL(2,\bC)$ of the ball $\frac{1}{2} D$ and map
this to $X$ by $g\rightarrow g(x_{0})$. The image obviously contains a neighbourhood
$N$ of $x_{0}$ and on $N$ we have $\phi\leq \frac{\oa-b}{2}+ f_{0}$. Then
Lemma 1 implies that
$b$ cannot be very large. In fact, if  the minimum of $f_{0}$ on $N$ is $2\pi
\ua$, we have $\phi\leq ( \frac{\oa}{2}- \ua)- \frac{b}{2}$ on $N$, so
$$ -M\leq \int_{N} \phi\  d\mu_{0} \leq ((\frac{\oa}{2}-\ua) - \frac{b}{2}) {\rm Vol}(N), $$
hence
$$ b\leq (\oa-2 \ua)+  \frac{2M}{{\rm Vol}(N)} $$
where $M$ is as in Lemma 1. This completes the proof of Theorem 1.

\section{The toric case}

Now suppose that $Z$ is a toric Fano manifold associated to a polytope $P\subset
V\cong \bR^{n}$. Thus the compact $n$-dimensional torus $T$ acts on $Z$. Suppose $\Gamma\subset SL(n,{\bf Z})$ is a finite group of symmetries
of $P$, with only one fixed point in $P$. For example  $Z$ could be 
the projective plane blown up in three general points, in which case $P$
is a hexagon and we can take $\Gamma$ to be a cyclic group of order $6$ (or
a dihedral group). Then there is a group  $\hat{T}$  which fits into a split exact
sequence
\begin{equation} 1 \rightarrow T \rightarrow \hat{T} \rightarrow \Gamma \rightarrow
1 \end{equation} and which  acts on $Z$. There is a unique $T$-orbit $M$ in $Z$ which is
preserved by the action of any lift of $\Gamma$ and we fix a basepoint $x_{0}$ in $M$.
\begin{prop}
The $\alpha$-invariant $\alpha_{G}(Z)$ is $1$.
\end{prop}
This is a theorem of Batyrev and Selinova \cite{kn:BS}. Song gave another proof, and established a converse,
in \cite{kn:S}. We will see that the argument can be made almost identical to that
of Theorem 1. We have a section $s$ of $K_{Z}^{-1}$ vanishing along a divisor
$D$ whose components  correspond to the $(n-1)$-dimensional faces of the
 boundary of $P$. For any $T$ invariant Hermitian metric on $K_{Z}^{-1}$
 the function $f=-\log \vert s\vert^{2}$ is $T$-invariant and, by considering
 the $T^{c}$-orbit
 of $x_{0}$, this induces a convex function $\tilde{f}$ on the dual space $V^{*}= T^{c}/T$. The polytope $P$ is the image of the derivative of  $\tilde{f}$.
Then $\Gamma$ acts on $V^{*}$ with a unique fixed point. For any $\hat{T}$-invariant
metric on $Z$, the function $\tilde{f}$ is preserved by the $\Gamma$-action
on $V^{*}$. From there the proof
proceeds exactly as before. The analogue of Lemma 1 holds with $\beta< 1$
since the local models for the zeros of $s$ are $f_{p}(z_{1},
\dots, z_{n})=0$ where $f_{p}(z_{1},\dots, z_{n})=z_{1}\dots z_{p}$ and 
$\vert f_{p}\vert^{-2\beta}$ is locally integrable for $\beta<1$.

\

Notice that, while the differential-geometric objects we work with are identical
in the two cases (convex functions on symmetric spaces $G^{c}/G$, invariant
under a finite group $\Gamma$ acting on $G^{c}/G$ with a single fixed point),
the formal structure of the group actions is different. In one case $\Gamma$
is a subgroup of $G$ and in the other it fits into an exact sequence (1).

\section{Review of the Mukai-Umemura manifold}

We give a brief summary of the most relevant parts of \cite{kn:MU}.

Start with a $2$-dimensional complex vector space $V$, so $\bP(V)=\bC\bP^{1}$
is the Riemann sphere. Fix a Hermitian form on $V$, so $SU(V)$ acts by isometries
on $\bP(V)$. 
Consider the configuration of vertices of a regular icosahedron as a  point $\delta$ in the symmetric
product ${\rm Sym}^{12}(\bP(V))$ and identify this with the projectivization of the
symmetric power $s^{12}=s^{12}(V)$ in the usual way. So $\Gamma\subset PSL(V)$ acts
on $\bP(s^{12}) $ fixing $\delta$. Now $-1\in SL(V)$ acts trivially on $s^{12}$,
so $PSL(V)$ acts on $s^{12}$ and $\Gamma$ fixes the $1$-dimensional subspace
of $s^{12}$ corresponding to $\delta$.  The action
on this $1$-dimensional space gives a character ${\Gamma}\rightarrow S^{1}$ and, since
$\Gamma$ is a perfect group this must be trivial. So if we take a
generator $\hat{\delta}$ for this line the stabiliser of $\hat{\delta}$ in
$PSL(V)$ is $\Gamma$. Thus the orbit  of $\hat{\delta}$ in the $13$-dimensional
complex vector space $s^{12}$ is a copy of $PSL(2,\bC)/\Gamma$. Now take
the closure of this orbit in the standard compactification $\bP^{13}= \bP(\bC\oplus
s^{12})$ of the complex vector space. This is a complex variety $X$ with
a $PSL(V)$ action induced by the linear action on $\bC\oplus s^{12}$ (trivial
on the first factor). It contains a dense open orbit, the orbit of $\hat{\delta}=x_{0}$. Mukai and
Umemura show that $X$ is smooth and that the embedding $X\rightarrow \bP^{13}$
is defined by  the complete linear system $-K_{X}$. So $X$ is Fano
and there is a $PSL(V)$-invariant section $s$ of $K_{X}^{-1}$ which vanishes
on the section at infinity $\Sigma\subset X$ (In other words, $H^{0}(K_{X}^{-1}) = \bC\oplus s^{12}$
as a representation of $PSL(V)$.)

Mukai and Umemura show that the infinity section $\Sigma=X\cap \bP(s^{12})$ can be described
as follows. The rational normal curve is the image of the map which assigns
to a point $\alpha\in \bP(V)$ the point $(\alpha,\alpha,\alpha,\dots,\alpha)$
in $\bP(s^{12})$, regarded as the symmetric product. In co-ordinates this
is the map which takes $\alpha$ to the co-efficients of the polynomial
$(z-\alpha)^{12}$. Now we map $\bP(V)\times \bP(V)$ to $\bP(s^{12})$ by mapping
a pair $(\alpha, \beta)$ to the point $(\alpha, \dots,\alpha, \beta)$ in
the symmetric product. In coordinates this is the map which takes $(\alpha,
\beta)$ to the co-efficients of the polynomial $(z-\alpha)^{11} (z-\beta)$.
This is obviously equivariant with respect to the $PSL(V)$ action on the
two spaces. The surface $\Sigma$ is the image of this map. The map is clearly
injective so, as a set, $\Sigma$ can be identified with $\bP(V)\times\bP(V)$.
Since $-c_{1}(K_{X})^{3}=22$, which has no cube factor, and  the divisor $\Sigma$ has
only one component, the section $s$ must vanish with multiplicity $1$ along
$\Sigma$. So what remains is to identify the form of the singularity of
$\Sigma$.

It is easy to see that the only singularities in $\Sigma$ occur along the
rational normal curve, which is the image of the diagonal in $\bP(V)\times
\bP(V)$. There are a number of ways of showing that the singularity at these
points has the given form. First, we can do a direct calculation.
 By
the $PSL(V)$ symmetry it suffices to look around any point on the diagonal,
so we can work in local co-ordinates with $\alpha, \beta$ close to $0\in
\bC$. Then, taking the co-efficients of $(z-\alpha)^{11}(z-\beta)$, the map is given in c-ordinates by
$$f(\alpha,\beta)= (c_{1} \alpha + \beta, c_{2} \alpha^{2} + c_{1} \beta\alpha,
c_{3}\alpha^{3} +c_{2} \alpha^{2} \beta, \dots);
$$
where the terms not written have order $4$ or more in $\alpha, \beta$ and
the $c_{i}$ are the binomial co-efficients. If we set $u= c_{1} \beta +\alpha$ the map can be written as
$$  f(u,\alpha)= (u, (c_{2}- c_{1}^{2})\alpha^{2} + c_{1} \alpha u, (c_{3}
- c_{1} c_{2}) \alpha^{3} + c_{2} \alpha^{2} u , \dots ). $$
The image of the curve $u=0$, which is a section transverse to the rational normal
curve, obviously has a cusp singularity. Using the fact that $PSL(V)$ acts
transitively on the rational normal curve it is not hard to deduce that the
surface has a singularity of the given form. 
 (Note  also that the surface $\Sigma$ is the surface
swept out by the lines tangent to the rational normal curve. For any projective
curve
$C$ with non-vanishing curvature and torsion  the surface swept out by the
tangent lines has a cusp singularity along $C$.)

For a second point of view, we consider the more general situation where
$PSL(2,\bC)$ acts on a $3$-dimensional complex manifold with a closed orbit
$C$ isomorphic to $\bP^{1}$. The stabiliser of a point in $C$ acts on the normal bundle of $C$
with a pair of weights $p,q$.
The Luna Slice Theorem gives a complete description of the possible $2$-dimensional
orbits, in a neighbourhood of the curve $C$. The closure of one of these orbits
is either smooth or modelled on the plane curve $z^{p}=w^{q}$, in a slice
transverse to $C$. So in the case at hand one can obtain the required result
by verifying that these weights are $2,3$.

For a third point of view we
just use
the facts that $\Sigma$ is an anti-canonical divisor in $X$ and that $\bP^{1}\times
\bP^{1}\rightarrow \Sigma$ is a bijective map which is an immersion away from
the diagonal. We can consider the  structure sheaf of $\Sigma$  as a sheaf ${\cal S}$ over the normalisation
$\bP^{1}\times \bP^{1}$. The quotient ${\cal O}_{\bP^{1}\times \bP^{1}}/{\cal S}$ is supported on the diagonal, As Alessio Corti explained to the author, the fact
that the dualising sheaf of $\Sigma$ is trivial forces the quotient
 to be the sheaf of sections of the co-normal bundle of
the diagonal in $\bP^{1}\times \bP^{1}$,  which implies in turn that $\Sigma$ has a cusp singularity. In other
words, the only way we can obtain an anti-canonical divisor in a Fano 3-fold
which is normalised by $\bP^{1}\times \bP^{1}$, in this fashion, is to have
a cusp singularity transverse to the diagonal.



\end{document}